\documentclass{amsart}
\usepackage[final]{epsfig}
\usepackage{graphics}
\usepackage{float}
\usepackage{amsmath}
\usepackage{amsfonts}
\usepackage{latexsym}
\usepackage{amssymb}
\usepackage{graphicx}
\usepackage{url}
\usepackage{epstopdf}
\usepackage{verbatim}
\usepackage[margin=1in]{geometry}
\usepackage{amsthm}
\usepackage{tikz}
\usepackage{caption}
\usepackage{tikz-cd}
\usepackage{enumitem}

\makeatletter
\newtheorem*{rep@theorem}{\rep@title}

\newcommand{\newreptheorem}[2]{%
\newenvironment{rep#1}[1]{%
 \def\rep@title{#2 \ref{##1}}%
 \begin{rep@theorem}}%
 {\end{rep@theorem}}}
\makeatother

\newtheorem{lemma}{Lemma}[section]
\newtheorem{proposition}[lemma]{Proposition}

\newtheorem{theorem}[lemma]{Theorem}
\newtheorem{definition}[lemma]{Definition}
\newtheorem{corollary}[lemma]{Corollary}

\newtheorem*{theorem*}{Theorem}
\newreptheorem{theorem}{Theorem}

\newcommand{\proofend}{$\Box$\bigskip}

\newcommand{\N}{{\mathbb N}}
\newcommand{\R}{{\mathbb R}}

\def\proof{\par{\it Proof}. \ignorespaces}

\newtheorem*{thm*}{Theorem}
\newtheorem*{thma*}{Theorem A}
\newtheorem*{thmb*}{Theorem B}
\newtheorem*{thmc*}{Theorem C}
\newtheorem*{thmd*}{Theorem D}

\begin{document}

\title{Determining the Maximum Difference between the Number of Atoms and number of Coatoms of a Bruhat interval of the Symmetric Group}

\author{E. Tsukerman}
\date{\today}

\begin{abstract}
We determine the largest difference between the number of atoms and number of coatoms of a Bruhat interval of $S_n$.
\end{abstract}

\address{Department of Mathematics, University of California,
Berkeley, CA 94720-3840}
\email{e.tsukerman@berkeley.edu}

\maketitle
\setcounter{tocdepth}{1}

\section{Introduction}

Much work has been done on understanding the structure of Bruhat intervals of the symmetric group (see, e.g., \cite{MR644668}, \cite{MR1970983} and \cite{BB} along with references therein).
Recently, particular interest has arisen in understanding the number of atoms and coatoms of Bruhat intervals of the symmetric group \cite{MR2221362, MR2727461}. There, the maximum number of atoms and coatoms of an interval of a given length is determined. In this note, we determine the largest difference between the number of atoms and the number of coatoms of a Bruhat interval of $S_n$.

Our main results are

\begin{thma*}

Let $\mathcal{I}$ be the set of intervals in $S_n$ and for $I \in \mathcal{I}$, let $a(I)$ and $c(I)$ denote the number of atoms and coatoms of $I$ respectively. Then
\[
\max_{I \in \mathcal{I}} c(I)-a(I) = \lfloor n^2/4 \rfloor -n+1.
\]
\end{thma*}

\begin{thmb*}
Let $n \geq 4$. An interval $I=[u,v] \subset S_n$ maximizes $c(I)-a(I)$ if and only if $c(I)=\lfloor n^2/4 \rfloor$ and $a(I)=n-1$.
\end{thmb*}

\section{Facts about Bruhat Intervals in $S_n$}

We will be needing the following definition and two results.

\begin{definition}\cite[Definition 4.9]{1406.5202}\label{defGat}
Let $u \leq v$ be permutations in $S_n$,
 and let $\overline{T}([u,v]):=\{t  \in T : u \lessdot ut \leq v\}$ 
and $\underline{T}([u,v]):=\{t  \in T : v \gtrdot vt \geq u\}$
be the 
transpositions labeling the cover relations corresponding
to the atoms and coatoms in the interval. Define a labeled graph $G^{at}$ (resp. $G^{coat}$) on $[n]$ such that $G^{at}$ (resp. $G^{coat}$) has an edge
between $a$ and $b$ if and only if 
the transposition $(a b) \in \overline{T}([u,v])$ (resp. $(a b) \in \underline{T}([u,v])$). Let $B_{u,v}^{at}$
be the partition of $[n]$ whose
blocks  are the connected components of $G^{at}$.
Similarly, define partition $B_{u,v}^{coat}$ whose
blocks  are the connected components of $G^{coat}$.
\end{definition}

The next result allows us to relate the atoms and coatoms of a Bruhat interval.

\begin{proposition}\cite[Proposition 4.10]{1406.5202} \label{atomPartition}
Let $[u,v] \subset S_n$.
The labeled graphs  $G^{at}$ and $G^{coat}$ have the same connected components.
\end{proposition}

The following result gives a sharp upper bound on the number of coatoms an interval of $S_n$ can have.

\begin{proposition}\label{maxAtom}\cite[Proposition 2.1]{MR2221362} 
For every positive integer $n$,
\[
\max_{v \in S_n} \# \underline{T}([1,v])=\lfloor n^2/4 \rfloor.
\]
\end{proposition}

The final result describes the permutations $v$ for which the number of coatoms is maximal.

\begin{proposition}\label{whenMaxAtom}\cite[Proposition 2.9]{MR2221362} 
For every positive integer $n$,
\[
\# \{v \in S_n \, \vert \, \#\underline{T}([1,v])=\lfloor n^2/4 \rfloor\}=\begin{cases}
n, & \text{ if } n \text{ is odd}; \\
n/2, & \text{ if } n \text{ is even}.
\end{cases}.
\]
Each such permutation has the form
\begin{align}\label{optTop}
v=[t+m+1,t+m+2,\ldots,n,t+1,t+2,\ldots,t+m,1,2,\ldots,t],
\end{align}
where $m \in \{ \lfloor n/2 \rfloor, \lceil n/2 \rceil\}$ and $1 \leq t \leq n-m$.
\end{proposition}

\section{Largest gap between number of atoms and coatoms of an interval in the Symmetric Group}

In this section, we consider the question of how large a gap can there be between the number of atoms and coatoms  of a Bruhat interval of the Symmetric group $S_n$. The first result is a simple inequality that will be used later in finding a maximum.

\begin{lemma} \label{floorLemma}
For all $k_1,k_2 \in \N$ with $k_i \geq 2$,
\[
\lfloor k_1^2/4 \rfloor+\lfloor k_2^2/4 \rfloor+1 < \lfloor (k_1+k_2)^2/4 \rfloor.
\]
\end{lemma}

\proof
We have
\[
\lfloor k_1^2/4 \rfloor+\lfloor k_2^2/4 \rfloor \leq \lfloor k_1^2/4+k_2^2/4 \rfloor.
\]
Therefore it suffices to observe that for $k_1,k_2 \geq 2$,
\[
k_1^2/4+k_2^2/4+1 < (k_1+k_2)^2/4=k_1^2/4+k_2^2/4+\frac{k_1 k_2}{2}.
\]
\proofend

We now prove the main result of this note, which states that the largest difference between the number of coatoms and atoms of an interval of $S_n$ is equal to $\lfloor n^2/4 \rfloor -n+1$.

\begin{theorem}\label{maxDiff}
Let $\mathcal{I}$ be the set of intervals in $S_n$ and for $I \in \mathcal{I}$, let $a(I)$ and $c(I)$ denote the number of atoms and coatoms of $I$ respectively. Then
\[
\max_{I \in \mathcal{I}} c(I)-a(I) = \lfloor n^2/4 \rfloor -n+1.
\]
\end{theorem}

\proof
Let $I \in \mathcal{I}$ and consider $G^{at}$ and $G^{coat}$ as in definition \ref{defGat}. 
By Proposition \ref{atomPartition}, $G^{at}$ and $G^{coat}$ have the same connected components. Let $\mathcal{K}_i$, $i=1,2,\ldots,m$, be the connected components of $G^{at}$ and $G^{coat}$, and let $k_i \geq 1$ denote their respective number of vertices. Let $p$ be the number of active components and $q=m-p$. By Proposition \ref{maxAtom},
\[
c(I) \leq \sum_{i=1}^m \lfloor k_i^2/4 \rfloor.
\]
Also,
\[
a(I) \geq \sum_{i=1}^m (k_i-1).
\]
Therefore
\begin{align} \label{eq1}
 c(I)-a(I) \leq \sum_{i=1}^m \lfloor k_i^2/4 \rfloor -(k_i-1).
\end{align}
We maximize (\ref{eq1}) over possible $\mathcal{K}_i$. Let $f(x)=\lfloor x^2/4 \rfloor-x+1$, so that
\[
c(I)-a(I) \leq \sum_{i=1}^m f(k_i).
\]
By Lemma \ref{floorLemma}, if $k_1,k_2 \geq 2$, then
\[
f(k_1+k_2)=f(k_1+k_2)+f(1) > f(k_1)+f(k_2).
\]
Therefore 
\[
\sum_{i=1}^m f(k_i) \leq qf(1)+(p-1)f(1)+f(n-q)=f(n-q).
\]
Since 
\[
\Delta[f](n)=f(n+1)-f(n)=\begin{cases}
\frac{n}{2}-1 & \text{if } n \text{ is even} \\
\frac{n+1}{2}-1 & \text{if } n \text{ is odd},
\end{cases}
\]
the function $f:\N \rightarrow \R$ is monotonically increasing. It follows that for every $I \in \mathcal{I}$,
\[
c(I)-a(I) \leq f(n).
\]
Next we show that the value $f(n)$ is attained for some interval $I=[u,v]$. We consider any $v$ of the form (\ref{optTop}). By Proposition \ref{whenMaxAtom}, the interval $[1,v]$ has $\lfloor n^2/4 \rfloor$ coatoms. The identity permutation has exactly $n-1$ elements covering it.
\proofend

\begin{theorem}
Let $n \geq 4$. An interval $I=[u,v] \subset S_n$ maximizes $c(I)-a(I)$ if and only if $c(I)=\lfloor n^2/4 \rfloor$ and $a(I)=n-1$.
\end{theorem}
\proof
From the proof of Theorem \ref{maxDiff}, 
\[
c(I)-a(I) \leq f(n-q).
\]
The assumption that $n \geq 4$ implies that $f(n)>f(n-q)$ for every $q>0$. Moreover, we know that the maximum value of $f(n)$ is attainable. Therefore $c(I)-a(I)$ is maximized only if $q=0$. So assume that $q=0$. Let $K$ be the single connected component of $G^{at}$ and $G^{coat}$ which contains $n$ vertices. Then
\[
c(I) \leq \lfloor n^2/4 \rfloor
\]
and
\[
a(I) \geq n-1.
\]
\proofend

\begin{corollary}
Let $n \geq 4$. Suppose that $I=[u,v] \subset S_n$ is an interval for which $c(I)-a(I)$ is maximized. Then $v$ is of the form (\ref{optTop}).
\end{corollary}

\proof
By Proposition \ref{whenMaxAtom}, the number of coatoms is $\lfloor n^2/4 \rfloor$ only for $v$ of the form (\ref{optTop}).
\proofend

A family of intervals for which the optimal value $c(I)-a(I) = \lfloor n^2/4 \rfloor -n+1$ is attained is given by
\[
I=[1,v]
\]
for $v$ as in (\ref{optTop}). There exist other intervals for which this maximum is attained. For example, in $S_4$, the intervals for which the maximum is attained are 
\[
[1234,3412],[1234,4231],[1243,4231],[2134,4231].
\]

\bigskip
{\bf Acknowledgments}.  
This material is based upon work supported by the National Science Foundation Graduate Research Fellowship under Grant No. DGE 1106400. Any opinion, findings, and conclusions or recommendations expressed in this material are those of the authors(s) and do not necessarily reflect the views of the National Science Foundation.

\bibliographystyle{alpha}
\bibliography{bibliography}



\end{document}